%% file: main.tex
\documentclass[letterpaper, 10 pt, conference]{ieeeconf}  % Comment this line out
\IEEEoverridecommandlockouts                              % This command is only
\input{avz_defs}

\usepackage{algorithmic}
\usepackage{epstopdf}
\usepackage{float}
\title{\LARGE \bf
Optimal Control for Scheduling and Pricing Intra-day Natural Gas Transport on Pipeline Networks
}

\author{Anatoly Zlotnik$^\dagger$, Kaarthik Sundar$^\star$, Aleksandr M. Rudkevich$^\ddag$, Aleksandr Beylin$^\ddag$, and Xindi Li$^\ast$
\thanks{$^\dag$azlotnik@lanl.gov, \,\, Applied Mathematics \& Plasma Physics, Los Alamos National Laboratory, Los Alamos, NM 87545}
\thanks{$^\star$kaarthik@lanl.gov \,\, Information Systems \& Modeling, Los Alamos National Laboratory, Los Alamos, NM 87545}
\thanks{$^\ddag$arudkevich@negll.com, abeylin@negll.com, \,\ Newton Energy Group, Boston, MA 02116}
\thanks{$^\ast$xindi@tcr-us.com \,\ Tabors Caramanis Rudkevich, Boston, MA 02116}}

\begin{document}

\maketitle
\thispagestyle{empty}
\pagestyle{empty}

%remove this!
%%%%%%%%%%%%%%%%%%%%%%%%%%%%%%%%%%%%%%%%%%%%%%%%%%%%%%%%%%%%%%%%%%%%%%%%%%%%%%%%
\begin{abstract}
We formulate an economic optimal control problem for transport of natural gas over a large-scale transmission pipeline network under transient flow conditions.  The objective is to maximize economic welfare for users of the pipeline system, who provide time-dependent price and quantity bids to purchase or supply gas at metered locations on a system with time-varying injections, withdrawals, and control actions of compressors and regulators.  Our formulation ensures that pipeline hydraulic limitations, compressor station constraints, operational factors, and pre-existing contracts for gas transport are satisfied.  A pipeline is modeled as a metric graph with gas dynamics partial differential equations on edges and coupling conditions at the nodes. These dynamic constraints are reduced using lumped elements to a sparse nonlinear differential algebraic equation system.  A highly efficient temporal discretization scheme for time-periodic formulations is introduced, which we extend to develop a rolling-horizon model-predictive control scheme.  We apply the  computational methodology to a pipeline system test network case study. In addition to the physical flow and compressor control solution, the optimization yields dual functions that we interpret as the time-dependent economic values of gas at each location in the network.  
\end{abstract}

%%%%%%%%%%%%%%%%%%%%%%%%%%%%%%%%%%%%%%%%%%%%%%%%%%%%%%%%%%%%%%%%%%%%%%%%%%%%%%%%
\section{Introduction} \label{secintro}

As electric power systems in many parts of the world increasingly rely on gas-fired generation, mechanisms for economically and operationally efficient coordination between the wholesale natural gas and electricity markets are of increasing interest \cite{unsihuay07a,liu11,zlotnik17tpwrs}.  System operators in both these sectors desire more efficient and reliable decision support tools that provide accurate price signals to inform operating and investment decisions \cite{hibbard12,tabors14}.  Power system operation and wholesale electricity pricing is currently conducted in organized optimization-based electricity markets administered by regional transmission organizations, so that prevalent electric energy prices are consistent with the physical capacity of the power grid \cite{pjm15}.  An optimization-based approach for scheduling natural gas flows throughout pipeline systems could enable computation of location- and time-dependent prices of natural gas that account for pipeline engineering factors, operational constraints, and the physics of gas flow \cite{zlotnik17psig}.  Efficient coordination between the two sectors could then be facilitated by the exchange of physical flow and price time-series between participants in the corresponding markets, in which prices are computed to be consistent with the physics of energy flow \cite{rudkevich17hicss}. The communicated physical data would be forecast or desired hourly energy consumption schedules, and pricing data would be bids and offers that reflect the willingness to transact payments for energy.

Optimization-based markets for physical flow scheduling and formation of location- and time- dependent pricing of natural gas are intended to address the needs of gas-fired generators, which may quickly change their fuel consumption \cite{zhao18shadow}.  Such a mechanism would therefore require accurate representation of hydraulic transients in gas pipelines within an optimization formulation.  It is well understood that compressibility of natural gas significantly affects the propagation of changes in pressures and flows throughout a large pipeline system, and therefore steady-state models, or sequences of such models, are insufficient to capture the effect of changes in mass within a section of pipe, or so-called ``line-pack'' \cite{carter03,liu11}.  Transient optimization, which refers to optimal control of gas pipeline dynamics, and pipeline model predictive control (MPC) have been proposed in a number of studies \cite{steinbach07pde,abbaspour07,rachford09,gopalakrishnan13mpc,zlotnik15cdc}, and there has been a resurgence of interest in recent years \cite{gugat18mip}.

In the previous transient optimization studies, the applicability to general network structures, scalability of the computational methods, and accuracy of models and solutions have presented challenges.  Recent work by the authors and collaborators has led to accurate and validated reduced pipeline dynamics partial differential equation (PDE) models, differential algebraic equation (DAE) discretization schemes, and problem formulations that are suitable for tractable, rapid pipeline transient optimization.  In particular, these recent studies have resulted in modeling concepts and dynamic system representations for general large-scale pipeline systems \cite{zlotnik15dscc,dyachenko17ss}, optimal control formulations \cite{zlotnik15cdc}, comparisons of various discretization schemes \cite{mak16acc,mak19ijoc}, extension to non-ideal gas modeling \cite{gyrya19staggered}, and validation of these models with respect to real data and commercial solvers \cite{zlotnik17psig}.  Although notable computational goals of accuracy, computational speed, and scalability have been achieved with these simulation and optimal control studies, operators of gas pipelines require decision support systems that can be perameterized with data collected from pipeline instruments, and solved on commodity computing platforms to provide information that can be used to improve their business processes. 

In this paper, we formulate an optimal control problem (OCP) for clearing an intra-day pipeline market using day-ahead, hourly physical flow and financial bids, whose solution provides an optimal flow schedule and hourly locational trade values (LTVs) of natural gas, while ensuring that pipeline hydraulic limitations, compressor station constraints, operational factors, and pre-existing shipping contracts are satisfied.  The formulation is intended to represent a secondary auction market for trading hourly deviations with respect to baseline (usually constant) flows that are agreed upon in a primary market.  We then describe reduction of the OCP to a nonlinear program (NLP) optimization formulation using previously developed model reduction techniques to perform a spatial discretization, and describe a novel, efficient, and well-conditioned time-discretization technique that implicitly encodes time-periodic boundary conditions.  Extension of this time-periodic formulation to non-periodic boundary conditions is proposed as a method to formulate economic model-predictive control for a pipeline market to be re-solved hourly in a rolling-horizon manner using look-ahead inputs over 24 to 72 hours, to yield hourly prices (taken as the first hour of a solution on the time horizon) \cite{rudkevich19hicss}.

The rest of the manuscript is organized as follows.  In Section \ref{sec:model}, we review modeling and model reduction of natural gas flow in pipeline networks with compressors.  Section \ref{sec:ocp} contains an OCP formulation for maximizing economic welfare for pipeline market participants subject to time-varying injections, withdrawals, and actions of compressors and regulators.  In Section \ref{sec:numeric}, we describe a collocation scheme for time-discretization of time-periodic OCPs that has sparsity and conditioning advantages over previously applied pseudospectral schemes \cite{mak19ijoc}. We then extend this time-periodic formulation to one that is suitable for use with non-time-periodic boundary data, and discuss implementation details. Section \ref{sec:example} describes computational results for an economic OCP case study for a test pipeline network, and Section \ref{sec:conc} contains a discussion of the results.

\section{Modeling of Gas Pipeline Network Dynamics} \label{sec:model}

In this section we review the standardized modeling of large-scale gas transmission pipelines that has proven tractable for transient optimization \cite{zlotnik15cdc,zlotnik17tpwrs,sundar18tcst}.  Compressible gas flow in a horizontal pipe with slow transients that do not cause waves or shocks can be described using a simplification of the one-dimensional Euler equations \cite{Thorley1987},
\begin{flalign}
& \partial_t \rho + \partial_x \varphi = 0 \quad \text{and} \quad a^2 \partial_x \rho = -\frac{\lambda}{2D} \frac{\varphi |\varphi|}{\rho}. & \label{eq:pde_1}
\end{flalign}
The right and left equations above capture conservation of mass and momentum, respectively. The variables $\rho$ and $\varphi$ are instantaneous gas density and mass flux, respectively, and are defined on the domain $[0,L] \times [0,T]$ where $L$ is the pipe length and $T$ is a time horizon. The term on the right hand side of the second equation aggregates friction effects, where the parameters are the Darcy-Wiesbach friction factor $\lambda$ and pipe diameter $D$.  We assume that gas pressure $p$ and density $\rho$ satisfy the ideal equation of state $p = a^2 \rho$ with $a^2 = ZR\boldsymbol{T}$, where $a$, $Z$, $R$, and $\boldsymbol{T}$, are the speed of sound, gas compressibility factor, ideal gas constant, and constant temperature, respectively.  Multiple studies have supported the use of this simplification in the regime of slow transients \cite{osiadacz84,herty10}.  We use the ideal gas approximation for simplicity of exposition, though extension to non-ideal gas modeling is straightforward \cite{gyrya19staggered}. Equation \eqref{eq:pde_1} has a unique solution when the initial and boundary conditions, i.e., one of $\rho(0,t) = \ubar{\rho}(t)$ or $\varphi(0,t) = \ubar{\varphi}(t)$ and one of $\rho(L,t) = \bar{\rho}(t)$ or $\varphi(L,t) = \bar{\varphi}(t)$, are specified. For convenience and numerical conditioning, we apply the dimensional transformations
\begin{flalign} \label{eq:nondim}
&\hat{t}=\frac{t}{\ell_0/a}, \quad \hat{x}=\frac{x}{\ell_0}, \quad \hat{\rho}=\frac{\rho}{\rho_0}, \quad \hat{\varphi}=\frac{\varphi}{a\rho_0},&
\end{flalign}
where $\ell_0$ and $\rho_0$ are nominal length and density, to yield the non-dimensional equations
\begin{flalign}
& \partial_t \rho + \partial_x \varphi = 0 \quad \text{and} \quad \partial_x \rho = -\frac{\lambda \ell_0}{2D} \frac{\varphi |\varphi|}{\rho}. & \label{eq:pde_nondim}
\end{flalign}
The hat symbols above and henceforth are omitted.

Turbulent flow along a pipe creates friction that causes pressure to gradually decrease in the flow direction, so gas compressors must be used to maintain pressure and flow through the system.  We model compressor stations as controllers that change the density between station outlet and inlet, as a multiplicative ratio at a point $x = c$ with conservation of flow.  This is represented as $\rho(c^+, t) = \alpha(t) \cdot \rho(c^-,t)$ and $\varphi(c^+,t) = \varphi(c^-,t)$ where $\alpha(t)$ denotes the time-dependent compression ratio between suction (intake) and discharge (outlet) pressure.  

A large-scale gas transmission network can be modeled for transient analysis as a set of edges (representing pipes) that are connected at nodes (representing junctions) where the gas flow can be compressed, withdrawn from, or injected into the system.  We consider the system as a connected directed metric graph $(\mathcal V, \mathcal E)$ where $\mathcal V$ and $\mathcal E$ represent the sets of nodes and edges, respectively, where $(i,j) \in \mathcal E$ represents an edge that connects nodes $i, j\in \mathcal V$. The system state is given by $\rho_{ij}$ and $\varphi_{ij}$ for all $(i,j)\in\cE$, which denote the instantaneous density and per-area mass flux, respectively, on edge $(i,j) \in \mathcal E$ defined on the domain $[0, L_{ij}] \times [0,T]$.  Each edge $(i,j)$ is characterized by its length $L_{ij}$, diameter $D_{ij}$, and friction factor $\lambda_{ij}$, which constitute the metric. The cross-sectional area of a pipe is denoted by $X_{ij}$. For each edge $(i,j)$, the evolution of $\rho_{ij}$ and $\varphi_{ij}$ is given by \eqref{eq:pde_nondim}, i.e.,
\begin{flalign}
& \partial_t \rho_{ij} + \partial_x \varphi_{ij} = 0 \,\, \text{and} \,\, \partial_x \rho_{ij} = -\frac{\lambda_{ij} \ell_0}{2D_{ij}} \frac{\varphi_{ij} |\varphi_{ij}|}{\rho_{ij}} & \label{eq:pde_nondim_edge_both} 
\end{flalign}
Here the sign of $\varphi_{ij}$ indicates flow direction, and we may write $\varphi_{ij}(x_{ij},t)=-\varphi_{ji}(L_{ij}-x_{ij},t)$.  Each junction $i \in \mathcal V$ is associated with a time-dependent nodal density $\varrho_i(t): [0,T] \to \mathbb R_+$. The set of controllers is $\mathcal C\subset \mathcal E\times\{+,-\}$, where $(i,j)\equiv(i,j,+)\in\mathcal C$ is a controller located at node $i\in\mathcal V$ that adjusts density of gas flowing into edge $(i,j)\in\mathcal E$ in the $i \to j$ direction, while $(j,i)\equiv(i,j,-)\in\mathcal C$ denotes a controller located at node $j\in\mathcal V$ that adjusts density into edge $(i,j)\in\mathcal E$ in the direction $j\to i$. Compression is modeled as a multiplicative ratio $\ubar{\alpha}_{ij}:[0,T]\to\mathbb R_+$ for $\forall (i,j, +)\in\mathcal C$ and $\bar{\alpha}_{ij}:[0,T]\to\mathbb R_+$ for $\forall (i,j, -)\in\mathcal C$.

Let $\mathcal V_\sigma \subset \mathcal V$ denote the set of nodes where time-varying density is defined as $\sigma_j(t)$ at junction $j \in \mathcal V_\sigma$ and mass flow into the system is free. Mass flow withdrawals at the other junctions $j \in \mathcal V_q = \mathcal V \setminus \mathcal V_\sigma$ are denoted by $q_j(t)$. Borrowing from power systems nomenclature, we refer to the $\mathcal V_\sigma$ and $\mathcal V_q$ as the set of ``slack'' and ``non-slack'' nodes, respectively. 

\begin{figure}[t!]
\centering
\includegraphics[scale=1]{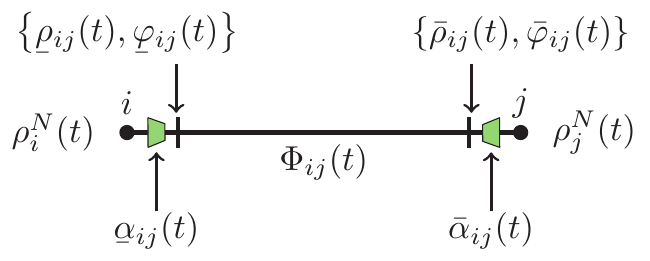}
\caption{The figure shows the densities and flows at the boundaries of each edge and the compression that can be applied at both the nodes $i$ and $j$.}
\label{fig:schematic}
\end{figure}

Nodal balance equations characterize the boundary conditions for the dynamics in Eq. \eqref{eq:pde_nondim_edge_both}. We define densities and flows at edge domain boundaries by
\begin{subequations}
\begin{flalign}
& \ubar{\rho}_{ij}(t) \triangleq \rho_{ij}(t, 0), \quad \bar{\rho}_{ij}(t) \triangleq \rho_{ij}(t, L_{ij}), & \label{eq:rhobar} \\
& \ubar{\varphi}_{ij}(t) \triangleq \varphi_{ij}(t, 0), \quad \bar{\varphi}_{ij}(t) \triangleq \varphi_{ij}(t, L_{ij}), & \label{eq:phibar}\\
& \text{and the nominal average edge flow as} \nonumber \\
& \Phi_{ij}(t) \triangleq \tfrac 12 (\ubar{\varphi}_{ij}(t) + \bar{\varphi}_{ij}(t)). &\label{eq:avg_phi}
\end{flalign}
\label{eq:bar_defn}
\end{subequations}
The above definitions are illustrated in Fig. \ref{fig:schematic} for a pipe joining two nodes $i$ and $j$.
Nodal balance laws are given using time-dependent compressor ratios $\ubar{\alpha}_{ij}(t)$ and $\bar{\alpha}_{ij}(t)$ , gas withdrawals $q_j(t)$, and supply densities $\sigma_j(t)$ as 
\begin{subequations}
\begin{flalign}
& \ubar{\rho}_{ij}(t) = \ubar{\alpha}_{ij}(t)\varrho_i(t), \, \forall \, (i,j) \in \mathcal E, & \label{eq:nodal_density_balance_1}\\
& \bar{\rho}_{ij}(t) = \bar{\alpha}_{ij}(t)\varrho_j(t), \, \forall \, (i,j) \in \mathcal E, & \label{eq:nodal_density_balance_2}\\
& q_j(t) =\sum_{i\in\mathcal V_q}X_{ij} \bar{\varphi}_{ij}(t)- \sum_{k\in\mathcal V_q}X_{jk}\ubar{\varphi}_{jk}(t), \,\forall\, j\in\mathcal V_q, & \label{eq:flow_balance} \\
& \ubar{\rho}_{ij}(t)=\sigma_i(t), \,\forall\, i\in\mathcal V_\sigma. & \label{eq:slack_pressure}
\end{flalign}
\label{eq:nodal_balance}
\end{subequations}
The optimized functions are time-varying compressor ratios $\{\ubar{\alpha}_{ij}, \bar{\alpha}_{ij}\}_{(i,j)\in\mathcal C}$ and the nodal gas withdrawals $\{q_j\}_{j\in\mathcal V_q}$.  Time-varying pressures at the slack nodes are given. We may omit dependence on time for ease of exposition.

The dynamics of gas flow over a pipeline network can be approximated by a reduced order nodal dynamics model. A lumped element approximation is made for \eqref{eq:pde_nondim_edge_both} on each edge, and equations \eqref{eq:nodal_balance} are written in terms of nodal density $\varrho_i$ for every $i \in \mathcal V$. This reduction extends the previous modeling work  \cite{grundel13a,zlotnik15cdc,zlotnik16ecc}. We consider a \textit{refinement} of a directed graph $(\mathcal V, \mathcal E)$ with edge lengths $L_{ij}$ to be $ (\hat{\mathcal V}, \hat{\mathcal E})$ when edges $(i,j)\in \hat{\mathcal E}$ are constructed by adding extra nodes to subdivide the edges of $\mathcal E$ such that the length $\hat L_{ij}$ of a new edge $(i,j) \in \hat{\mathcal E}$  satisfies $\hat L_{ij} < \Delta$.  The reduced model is shown to accurately resolve the PDE dynamics on a pipe when $\Delta$ is sufficiently small \cite{zlotnik15dscc}.  Lumping dynamics \eqref{eq:pde_nondim_edge_both} for each pipe segment $(i,j) \in \hat{\mathcal E}$ in the refined graph yields
\begin{subequations}
\begin{flalign}
& \int_0^L(\partial_t\rho_{ij}+\partial_x\varphi_{ij}) \,dx=0, & \label{eq:mass_int} \\
& \int_0^L(\partial_x \rho_{ij})\,dx = -\frac{\lambda_{ij}\ell_0}{2D_{ij}}\int_0^L \frac{\varphi_{ij} |\varphi_{ij}|}{\rho_{ij}}\,dx. & \label{eq:momentum_int}
\end{flalign}
\label{eq:pde_int}
\end{subequations}
The above integrals of $\partial_t$, $\partial_x$, and nonlinear terms are evaluated using the trapezoid rule, the fundamental theorem of calculus, and averaging variables, respectively, yielding
\begin{subequations}
\begin{flalign}
& \frac{L}{2}(\dot{\ubar{\rho}}_{ij}+\dot{\bar{\rho}}_{ij}) = \ubar{\varphi}_{ij}-\bar{\varphi}_{ij}, & \label{eq:mass_disc}\\
& \ubar{\rho}_{ij}-\bar{\rho}_{ij} = -\frac{\lambda_{ij}\ell_0 L}{4D_{ij}} \frac{(\ubar{\varphi}_{ij}+\bar{\varphi}_{ij}) |\ubar{\varphi}_{ij}+\bar{\varphi}_{ij}|}{\ubar{\rho}_{ij}+\bar{\rho}_{ij}}. & \label{eq:momentum_disc}
\end{flalign}
\label{eq:ode_disc}
\end{subequations}
%We remark that Eq. \eqref{eq:momentum_disc} can equivalently be written using the dissipation function as 
% \begin{flalign}
% & \Phi_{ij} + f_{\mu(ij)}\left(t, \frac{\ubar{\rho}_{ij} + \bar{\rho}_{ij}}{2}, \frac{\ubar{\rho}_{ij} - \bar{\rho}_{ij}}L\right) = 0 & \label{eq:momentum_f}
% \end{flalign}
% where $\Phi_{ij} = \frac 12 (\ubar{\varphi}_{ij} + \bar{\varphi}_{ij})$.
The equations \eqref{eq:ode_disc} and nodal balance laws \eqref{eq:nodal_balance} then reduce to DAE system:
\begin{subequations}
\begin{flalign}
& \frac{L}{2}(\dot{\ubar{\rho}}_{ij}+\dot{\bar{\rho}}_{ij}) = \ubar{\varphi}_{ij}-\bar{\varphi}_{ij},\, \forall \, (i,j) \in \hat{\mathcal E}& \label{eq:dae0c}\\
& \ubar{\rho}_{ij}-\bar{\rho}_{ij} = -\frac{\lambda_{ij}\ell_0 L}{D_{ij}} \frac{\Phi_{ij}|\Phi_{ij}|}{(\ubar{\rho}_{ij}+\bar{\rho}_{ij})},  \, \forall \, (i,j) \in \hat{\mathcal E} & \label{eq:dae0d} \\
& \ubar{\rho}_{ij} = \ubar{\alpha}_{ij}\varrho_i, \,\bar{\rho}_{ij} = \bar{\alpha}_{ij}\varrho_i, \, \forall \, (i,j) \in \hat{\mathcal E}, &\label{eq:dae0a}\\
& q_j = \sum_{i\in\hat{\mathcal V}_q}X_{ij} \bar{\varphi}_{ij} - \sum_{k\in\hat{\mathcal V}_q}X_{jk}\ubar{\varphi}_{jk}, \, \forall \, j\in\hat{\mathcal V}_q,  \label{eq:dae0b} & \\
& \ubar{\rho}_{ij} = \sigma_i, \,\forall \, i\in\hat{\mathcal V}_\sigma. & \label{eq:dae0e}
\end{flalign}
\label{eq:dae}
\end{subequations}
Eq. \eqref{eq:dae0a} represents continuity of density at junctions with jumps in the case of compression or regulation, Eq. \eqref{eq:dae0b} represents flow balance at junctions, and Eqs. \eqref{eq:dae0c}-\eqref{eq:dae0d} represent flow dynamics on each segment.

The DAE system in \eqref{eq:dae} can be written in matrix-vector form as follows. We enumerate the set of nodes in the set $\hat{\mathcal V}$ according to a fixed ordering, where non-slack nodes $\hat{\mathcal V}_q$ are ordered after the slack nodes, $\hat{\mathcal V}_\sigma$. Each node in $\hat{\mathcal V}$ is assigned an index $[\hat{\mathcal V}] := \{1,\dots, |\hat{\mathcal V}|\}$ according to the ordering. Each edge is assigned an index in $[\hat{\mathcal E}] := \{1, \dots, |\hat{\mathcal E}|\}$ and we define the map $\pi_e:\hat{\mathcal E} \to [\hat{\mathcal E}]$ that maps each edge to this ordering. Bold font henceforth represents vectors.

We now let $\bm \varrho = (\varrho_1, \varrho_2, \dots, \varrho_{|\hat{\mathcal V}|})^{\intercal}$ denote the nodal density state vector. 
Equation \eqref{eq:dae0a} will be used to state \eqref{eq:dae0c}-\eqref{eq:dae0d} in terms of nodal densities $\bm \varrho$.  We then define state vectors $\ubar{\bm \varphi}=(\ubar{\varphi}_1,\ldots,\ubar{\varphi}_{|\hat{\mathcal E}|})^\intercal$ and $\bar{\bm \varphi}=(\bar{\varphi}_1,\ldots,\bar{\varphi}_{|\hat{\mathcal E}|})^\intercal$, where $\ubar{\varphi}_k$ and $\bar{\varphi}_k$ are indexed by $k=\pi_e(ij)$. We denote $\bm \Phi = \tfrac 12 (\ubar{\bm \varphi} + \bar{\bm \varphi})$ as the vector of average flows on edges.

We now define the incidence matrix of the full refined graph $(\hat{\mathcal V}, \hat{\mathcal E})$, acting
$A:\mathbb R^{|\hat{\mathcal E}|}\to\mathbb R^{|\hat{\mathcal V}|}$, by
\begin{flalign} \label{eq:incidence0}
&A_{ik} = \left\{ \begin{array}{ll}  1 & \text{edge $k=\pi_e(ij)$ enters node $i$,} \\ -1 & \text{edge $k=\pi_e(ij)$ leaves node $i$,} \\ 0 & \text{else} \end{array}\right. &
\end{flalign}
and a weighted incidence matrix $B:\mathbb R^{|\hat{\mathcal E}|}\to\mathbb R^{|\hat{\mathcal V}|}$ given by
\begin{flalign} \label{eq:incidence0_w}
&B_{ik} = \left\{ \begin{array}{ll}  \bar{\alpha}_{ij} & \text{edge $k=\pi_e(ij)$ enters node $i$,} \\ -\ubar{\alpha}_{ij} & \text{edge $k=\pi_e(ij)$ leaves node $i$,} \\ 0 & \text{else}, \end{array}\right. &
\end{flalign}
where $\operatorname{sign}(B)=A$.  Here the compressor controls are embedded within the matrix $B$. A vector of withdrawal fluxes is defined by $\bm q=(q_1,\ldots,q_M)^T$ with $M=|\hat{\mathcal V}_q|$, where $q_k$ is negative if an injection. We also define the slack node densities as $\bm \sigma=(\sigma_1,\ldots,\sigma_b)^{\intercal}=\{\varrho_j\}_{j\in\hat{\mathcal V}_\sigma}$, where $b=|\hat{\mathcal V}_\sigma|$, and non-slack (demand) node densities as $\bm \rho=(\rho_1,\ldots,\rho_M)^T=\{\varrho_j\}_{j\in\hat{\mathcal V}_q}$, so that $b+M=|\hat{\mathcal V}|$. Note that $\bm \sigma$, $\bm \rho$ and $\bm \varrho$ are related by $\bm \varrho = (\bm \sigma, \bm \rho)^\intercal$, because of the choice of node ordering $\hat{\mathcal V}$. We let $A_\sigma,B_\sigma\in\mathbb R^{b\times |\hat{\mathcal E}|}$ denote the sub-matrices of rows of $A$ and $B$ corresponding to $\hat{\mathcal V}_\sigma$, and let $A_q,B_q\in\mathbb R^{M\times |\hat{\mathcal E}|}$ similarly correspond to $\hat{\mathcal V}_q$.  
%Next, let $A_L$ and $A_0$ denote the positive and negative parts of $A_d$, so that $A_d=A_L+A_0$.  
We then define the diagonal matrices $\Lambda,K, \boldsymbol{X}\in\mathbb R^{|\hat{\mathcal E}|\times |\hat{\mathcal E}|}$ by $\Lambda_{kk}=L_k$,  $K_{kk}=\ell_0\lambda_k/D_k$, and $\boldsymbol{X}_{kk} = X_k$ where $L_k$, $\lambda_k$, $D_k$, and $A_k$ are the non-dimensional length, friction factor, diameter, and cross-sectional area of edge $k=\pi_e(ij)$. Using this notation, \eqref{eq:dae} can be rewritten as a DAE system:
\begin{subequations}
\begin{flalign}
& |A_q| \boldsymbol{X} \Lambda |B_q^\intercal|\dot{\bm \rho} = 4(A_q \boldsymbol{X} \bm \Phi - \bm q) - |A_q| X \Lambda |B_\sigma^\intercal| \dot{\bm \sigma},  & \label{eq:dae1a} \\
& \Lambda K \bm \Phi \odot \bm \Phi =  -B^\intercal \bm \varrho \odot |B^\intercal| \bm \varrho, & \label{eq:dae1b} 
\end{flalign}
\label{eq:dae_final}
\end{subequations}
where the operator $\odot$ represents the Hadamard product.   Here, the gas withdrawals are $\bm q \in \mathbb R^M$, slack node densities are $\bm \sigma \in \mathbb R_+^b$, compression ratios are $\ubar{\alpha}_{ij}, \bar{\alpha}_{ij} \in \mathcal C$, and $\bm \rho \in \mathbb R_+^M$ and $\bm \Phi \in \mathbb R^{|\hat{\mathcal E}|}$ denote the system state.

To derive equations \eqref{eq:dae_final}, we rewrite Eq. \eqref{eq:dae0b} in matrix form as $\bm q = \bar{A}_q \boldsymbol{X} \bar{\bm \varphi} + \ubar{A}_q \boldsymbol{X} \ubar{\bm \varphi}$
where $\bar{A}_q$ and $\ubar{A}_q$ are the positive and negative parts of the matrix $A_q$, respectively. We now define $\bm \Phi_{-} = \frac 12 (\bar{\bm \varphi} - \ubar{\bm \varphi})$.  The Eq. \eqref{eq:dae0b} can then be rewritten as in the transformed variables $\bm \Phi$ and $\bm \Phi_{-}$ as
\begin{flalign}
& \bm q = A_q \boldsymbol{X} \bm \Phi + |A_q| \boldsymbol{X} \bm \Phi_{-}. & \label{eq:d_veca}  
\end{flalign}
Equations \eqref{eq:dae0a}, \eqref{eq:dae0c}, and \eqref{eq:dae0e} together with the definition $\bm \Phi_{-}$ can be equivalently represented using the matrix equation
\begin{flalign}
& |B_\sigma^\intercal| \dot{\bm \sigma} + |B_q^\intercal| \dot{\bm \rho} = -4 \Lambda^{-1} \bm \Phi_{-}.  & \label{eq:rho_dynamics_vec}
\end{flalign}
Substituting Eq. \eqref{eq:d_veca} into  \eqref{eq:rho_dynamics_vec} and eliminating $\bm \Phi_{-}$ yields \eqref{eq:dae1a}. Eq. \eqref{eq:dae0d} can be rewritten as \begin{flalign}
& \ubar{\rho}_{ij}^2 - \bar{\rho}_{ij}^2 = -\frac{\lambda \ell_0 L}{D_{ij}} \Phi_{ij}|\Phi_{ij}|, \, \forall (i,j) \in \hat{\mathcal E}. & \label{eq:phi_vec}
\end{flalign}
With equation \eqref{eq:dae0c} and the definitions of $B$, $\Lambda$, and $K$,  equation \eqref{eq:phi_vec} can be  written in matrix form as \eqref{eq:dae1b}.

\section{Economic Optimal Control Problem} \label{sec:ocp}

We formulate an economic OCP where \eqref{eq:pde_nondim_edge_both} and \eqref{eq:nodal_balance} are the dynamic constraints, for which we henceforth use the nodal equations.  In addition, we require several inequalities that arise from engineering limitations on the pipeline system.  First, there is a maximum allowable operating pressure (MAOP) at each point in the system, expressed as $\rho_{ij}(t,x)  \leq \ubar{\rho}_{ij}^{\max}$ for $\fA x\in[0,L_{ij}]$ and $\fA (i,j)\!\in\!\cE$.  These constraints may be enforced only at the endpoints of each pipe, because friction effects of turbulent flow subject to slowly varying transients result in monotone decrease of pressure along the direction of flow \cite{zlotnik16ecc}.  Minimum  pressure must be maintained at nodes, per contractual agreement.  We express these constraints as
\begin{subequations}
\begin{align} 
\ubar{\rho}_{ij}(t), \bar{\rho}_{ij}(t)  &\leq \ubar{\rho}_{ij}^{\max}, \quad \fA (i,j)\!\in\!\cE \\
\rho_i(t) &\geq \rho_i^{\min},  \quad \fA i\!\in\!\cV. 
\end{align}
\label{eq:plim0}
\end{subequations}
Next, the energy (or power) used by compressors is constrained by the inequalities
\begin{align}
\underline{\varepsilon}_{ij}|\underline{\phi}_{ij}(t)| \bp{(\underline{\alpha}_{ij}(t))^h-1} \leq \underline{E}_{ij}^{\max}, \quad (\underline{i,j})\in\underline{\cC}, \label{eq:comppow1a}\\
\overline{\varepsilon}_{ij}|\overline{\phi}_{ij}(t)| \bp{(\overline{\alpha}_{ij}(t))^h-1} \leq \overline{E}_{ij}^{\max}, \quad (\overline{i,j})\in\overline{\cC}, \label{eq:comppow1b}
\end{align}
with $h=(\gamma-1)/\gamma<1$ and where $\underline{\varepsilon}_{ij}$ and $\overline{\varepsilon}_{ij}$ correspond to $\varepsilon=(286.76 \cdot T_1) / (e_a \cdot e_m \cdot G \cdot h)$ for $(\underline{i,j})$ and $(\overline{i,j})$, respectively, where $T_1$, $e_a$, $e_m$, and $G$ are the discharge temperature, adiabatic and mechanical efficiencies, and gas gravity, respectively \cite{menon05}. We assume that compressor stations are designed and operated only to boost pressure, so
\begin{align} 
\underline{\alpha}_{ij}(t) \geq 1, \quad \overline{\alpha}_{ij}(t) \geq 1, \quad \fA (i,j)\in\cE. \label{eq:compmin}
\end{align}

The OCP of interest represents a two-sided single auction market that maximizes total surplus over injection and withdrawal schedules.  We define market surplus as the sum of producer (supplier) surplus and consumer (buyer) surplus.  Producer surplus occurs when the price a producer receives exceeds the value that they are willing to accept for the goods they sell.  Similarly, consumer surplus occurs when the price the consumer pays for good is below the value they are willing to offer.  Market surplus is the sum of individual surpluses for all consumers and producers who participate in the market.  We formulate an objective function similar to that used in previous studies \cite{oneill87,read12}.  The focus is on optimizing flows and pricing the value of gas deliveries as a function of time over an optimization horizon $[0,T]$, where $T$ is on the order of 12 to 72 hours.

\begin{figure}[t]
\centering
\includegraphics[width=.3\linewidth]{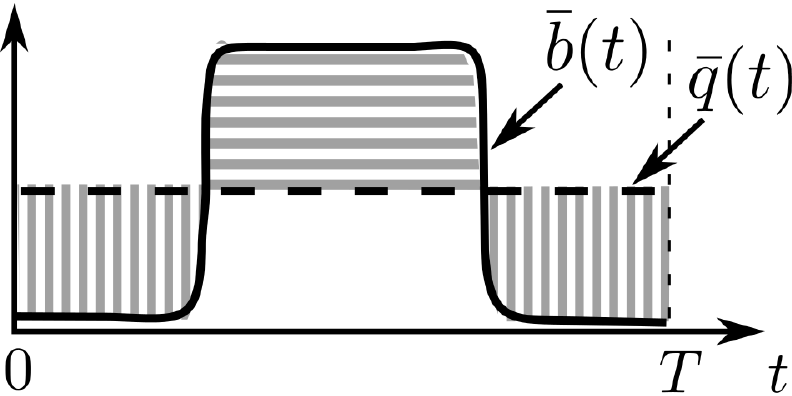} \,\,\,
\includegraphics[width=.3\linewidth]{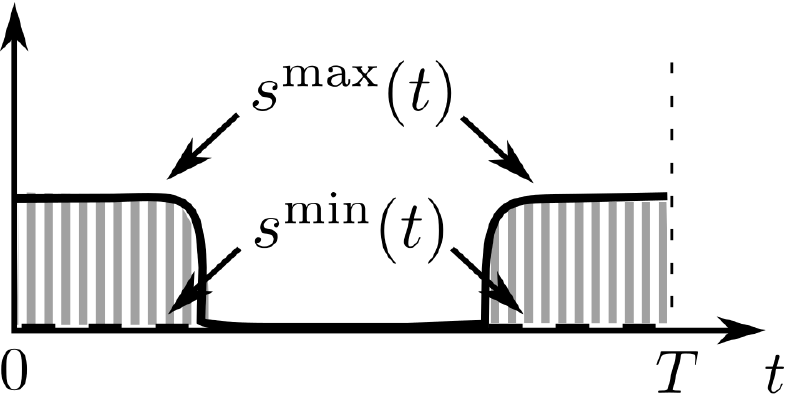} \,\,\,
\includegraphics[width=.3\linewidth]{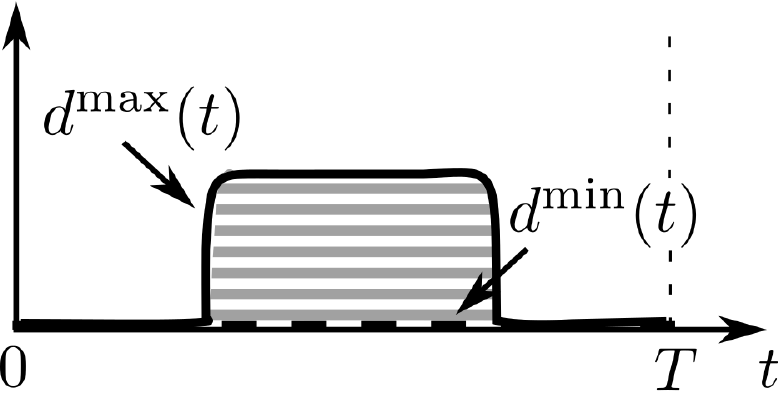}  \caption{Left: baseline flow $\bar{q}$ (dashed line) and quantity bid $\bar{b}$ (solid line), with surplus (vertical stripes) and deficit (vertical stripes) indicated; Center: bounds $d^{\min}$ and $d^{\max}$ on purchase bids for balancing, with feasible function space indicated by hatching; Right: bounds $s^{\min}$ and $s^{\max}$ on sell offers for balancing, with feasible function space indicated by hatching.} \label{fig:boundgen0}
\end{figure}

In order to account for the possibility of multiple customers and bidding structures at a single physical location, we introduce the set $\cG$ of transfer nodes in addition to the set of network nodes $\cV$.  The transfer nodes enumerate the notional receipt or delivery points associated with network nodes in $\cV$.  Each supplier is considered to be injecting gas at a unique transfer node $m\in\cG$, and each consumer withdraws gas at a unique node as well.  Each node $m\in\cG$ can represent only one supplier or consumer, and is associated with a unique network node $j(m)\in\cV$.  The set of transfer nodes connected to a node $j\in \cV$ is denoted by
\begin{align}
\partial_{g}j&=\left\{ m\in \cG\mid j(m)\in \cV\right\} \subset \cG.
\end{align}
We suppose that each slack node $j\in\cV_\sigma$ represents a single supplier transfer node where density $\sigma_j(t)$ is specified.  We use $\bar{q}_j(t)$ to denote the primary baseline flow withdrawal at node $j\in\cV$ about which a secondary auction is to take place.   The baseline profiles are assumed to have been agreed on based upon previously existing contracts, nominally for constant flow over the optimization period, and these withdrawals satisfy
\begin{align} \label{eq:injectbal1}
\sum_{j\in\cV} \int_0^T\bar{q}_j(t)\rd t = 0,
\end{align}
because they represent the outcome of a primary market mechanism.  These baseline withdrawal $\bar{q}_j$ for a physical node $j\in\cV$ can be decomposed into baseline supplies $\bar{s}_m$ and demands $\bar{d}_m$ at connected transfer nodes $m\in\partial_g j$, so that
\begin{align} \label{eq:inject1}
\bar{q}_j(t)&= \sum_{m\in\partial_g j} (\bar{d}_m(t)-\bar{s}_m(t)), \quad \fA \, t\in[0,T].
\end{align}
Unless baseline profiles are constant, their balancing is not necessarily instantaneous and must hold only as an integral over the planning horizon.  Then, the net (instantaneous) optimized variations in demand and supply at a transfer node $m$ with respect to baseline profiles $\bar{s}_m(t)$ and $\bar{d}_m(t)$ are represented by $d_m(t)$ and $s_m(t)$, respectively.  The total flow injection at a node is then formulated as
\begin{align} \label{eq:inject2}
q_j(t)&= \bar{q}_j(t) + \hat{d}_j(t)-\hat{s}_j(t), \quad \fA \, t\in[0,T],
\end{align}
where we denote
\begin{align} \label{eq:netimport0}
\hat{d}_j(t) = \sum_{m\in\partial_g j} d_m(t), \qquad \hat{s}_j(t) = \sum_{m\in\partial_g j} s_m(t).
\end{align}
Supplier limitations and consumer capacities at each $m\in\cG$ are subject to minimum and maximum constraints, which may depend on time, and are given by
\begin{subequations}
\begin{align} 
s_{m}^{\min}(t) \leq s_m(t) \leq s_{m}^{\max}(t), \quad \fA m\in\cG, \\ d_{m}^{\min}(t) \leq d_m(t) \leq d_{m}^{\max}(t), \quad \fA m\in\cG.
\end{align}
\label{eq:caplim0}
\end{subequations}
Because the optimized supplies $\hat{s}_j$ and deliveries $\hat{d}_j$ for each physical node $j\in\cV$ are specified with respect to the baseline flow $\bar{q}$, the bounds on the constraints in \eqref{eq:caplim0} are determined by the minimum and maximum deviations.  As an example, a procedure for generating the bound functions for these constraints given a baseline flow $\bar{q}$ and a desired quantity bid $\bar{b}$ of a single transfer node is illustrated in Figure \ref{fig:boundgen0}.
All suppliers and consumers place bids into the market that consist of minimum and maximum limits on supplies and consumptions, as well as offer prices $c_m^s(t)$ and bid prices $c_m^d(t)$, respectively, specified for each transfer node.  Thus, the parameters required to define inputs to the auction market are $c_m^s(t)$ and $c_m^d(t)$ with dimensions of price per unit mass; and $s_m^{\min}(t)$, $s_m^{\max}(t)$, $d_m^{\min}(t)$, $d_m^{\max}(t)$, and $\bar{q}_j(t)$  with dimensions of mass per unit time (mass flow). The market surplus objective function is given by
\begin{align} \label{eq:mswobj0}
J_{\rM{MS}} & \teQ  \sum_{m\in\cG} \int_0^T c_m^d(t) d_m(t) \rd t \nonumber \\ & \qquad - \sum_{m\in\cG} \int_0^T c_m^s(t) s_m(t) \rd t.
\end{align}
With the above collection of engineering and physical constraints, the optimal control formulation is
\begin{equation} \label{prob:msw0}
\begin{array}{ll}
\!\!\!\! \mathrm{max}  & J_{\rM{MS}} \teQ \text{max market surplus objective } \eqref{eq:mswobj0} \\
\!\!\!\! \text{s.t.} & \text{system dynamics } \eqref{eq:dae_final} \\
& \text{market integration } \eqref{eq:inject1}-\eqref{eq:netimport0} \\
& \text{compressor limits } \eqref{eq:comppow1a}-\eqref{eq:compmin} \\
& \text{pressure limits } \eqref{eq:plim0} \\
& \text{supply and demand limits } \eqref{eq:caplim0}
\end{array}
\end{equation}

Conceptually, after solving the OCP \eqref{prob:msw0} we wish to compute estimates of the value of natural gas at physical nodes throughout the system, as functions of time.  For the optimal solution, the values of the Lagrange multipliers that correspond to satisfaction of the equality constraint \eqref{eq:inject1} at a physical node $j\in\cV_q$, which we denote as $\lambda_j(t)$, represent sensitivity of the market surplus objective function value in \eqref{eq:mswobj0} to changes in nodal withdrawals $\hat{d}_j(t)$ or $\hat{s}_j(t)$.  The dual functions $\lambda_j(t)$ can then be interpreted as the incremental economic value of gas flow leaving the system from node $j$ at time $t$, when considering the market auction over the entire optimization horizon $T$.

For the time horizon $T$, we require that the state variables $\varphi_{ij}$ and $\rho_{ij}$ are time-periodic in order for the dynamic constraints to be well-posed, and time-periodicity also has to be imposed on the control and parameter functions $\{\ubar{\alpha}_{ij}, \bar{\alpha}_{ij}\}_{(i,j)\in\mathcal C}$, $\{q_j\}_{j\in\mathcal V_q}$, and $\{\sigma_j\}_{j\in\mathcal V_\sigma}$ as given by \eqref{eq:termcon1c}--\eqref{eq:termcon1e} (see \cite{zlotnik15cdc}).
This yields the terminal conditions 
\begin{subequations}
\begin{flalign}
&\rho_{ij}(0,x_{ij})=\rho_{ij}(T,x_{ij}),  \,\forall\, (i,j)\in\mathcal E, &\label{eq:termcon1a}\\
&\phi_{ij}(0,x_{ij})=\phi_{ij}(T,x_{ij}),  \,\forall\, (i,j)\in\mathcal E, &\label{eq:termcon1b}\\
&\ubar{\alpha}_{ij}(0)=\ubar{\alpha}_{ij}(T), \, \bar{\alpha}_{ij}(0)=\bar{\alpha}_{ij}(T),  \,\forall\, (i,j)\in\mathcal C, &\label{eq:termcon1c} \\
&q_j(0)= q_j(T), \,\forall\, j \in \mathcal V_q, &\label{eq:termcon1d} \\
&\sigma_j(0)=\sigma_j(T), \,\forall\, j \in \mathcal V_\sigma. &\label{eq:termcon1e}
\end{flalign}
\label{eq:termcon}
\end{subequations}
We formulate problem \eqref{prob:msw0} with time-periodic boundary conditions for conceptual and computational well-posedness. Conceptually, without some specification of the initial and terminal conditions, these states could be produced by the solver in unpredictable ways.  We address the computational details in the following section.

\section{Computational Approach} \label{sec:numeric}

A widely-used approach for transcribing OCPs to nonlinear programs involves pseudospectral approximation \cite{canuto06}, such as the Legendre-Gauss-Lobatto scheme \cite{ross03}, which we have applied to optimal control of gas pipeline networks in a previous study \cite{zlotnik15cdc}.  Here, we present a scheme specifically constructed for optimal control subject to a time-periodicity constraint, which uses a uniform collocation grid on the circular time domain.  With this time-periodic formulation and corresponding discretization, there is no issue of a Gibbs phenomenon that causes poorly conditioned approximation of OCP formulations near initial and terminal time points (see p. 44 of \cite{huntington07}).

We consider an OCP, or optimization problem in function space, of the form
\begin{subequations}
\begin{align}
	\min_u\ \ & J(x,u)=\int_0^T \mathcal{L}(t,x(t),u(t))dt, \label{eq:ocp0a} \\
	{\rm s.t.}\ \ & f(t,x(t),\dot{x}(t),u(t),h(t))=0, \label{eq:ocp0b} \\
	& g(x(t),u(t),h(t))\leq 0, \label{eq:ocp0c} \\
    & x(0)=x(T), \,\, u(0)=u(T), \label{eq:ocp0d} 
\end{align}
\label{eq:ocp0}
\end{subequations}
on $\cT=[0,T]$.  Here $\mathcal{L}\in C^\kappa$ is in the space $C^\kappa$ of continuous functions with $\kappa$ classical derivatives, and the dynamic constraints $f\in C_n^{\kappa-1}$ are in the space $C_n^{\kappa-1}$ of $n$-vector valued $C^{\kappa-1}$ functions, with respect to the state, $x(t)\in\bR^n$, it's derivative $\dot{x}(t)\in\bR^n$, the control input, $u(t)\in\bR^m$, and a family of parameter functions $h(t)\in\bR^r$.  We suppose that the latter are $C^\kappa$ and given as periodic on the domain $[0,T]$, i.e. $h(0)=h(T)$. The function $g$ specifies path (inequality) constraints, and the state and control solutions are constrained to be time-periodic.   The admissible set for controls $u$ includes the $C_m^{\kappa}$ functions on $\cT$.  This problem \eqref{eq:ocp0} is a time-periodic DAE reformulation of the problem examined in \cite{ruths11cdc}.  Given known control functions $u(t)$, we have proved in related work that the solution $x(t)$ for such a formulation will be unique \cite{sundar18tcst}.  Note that because the formulation \eqref{eq:ocp0} has no initial and terminal state constraints, other than time-periodicity, it can be used when the state is unknown or only partially observable. 

Here we introduce a simple, direct collocation procedure for constructing a finite-dimensional NLP that approximates the problem \eqref{eq:ocp0}.  We use a local, piecewise-linear scheme where we approximate a time-periodic function $y$ on the domain $[0,T)$ using a set of $N$ uniformly spaced collocation points $t_k=T(k-1)/N$ for $k=1,\ldots,N$, with values $\bar{y}_k=y(t_k)$.  For $t\in[t_k,t_{k+1})$, the approximation is
\begin{equation}
	\label{eq:Ix} y(t) \approx  \wh{y}_k^N(t) = \bar{y}_k+(\bar{y}_{k+1}-\bar{y}_k)\cdot \frac{N}{T}\cdot (t-t_k),
\end{equation}
where $t_{N+1}\equiv t_1$, $\bar{y}_{N+1}\equiv y_1$, etc.  This scheme satisfies $y(t_k) = \wh{y}_N(t_k) = \bar{y}_k$, so the physical meaning of the interpolating coefficients $\bar{y}_k$ are clearly the values of the function $y$ at uniform collocation points.  The collocation points are chosen to lie uniformly on $[0,T]$ including the endpoints, which we may map to the unit circle $[0,2\pi)$ in which case the terminal time point is equivalent to the initial point. We evaluate the integral in \eqref{eq:ocp0a} and the derivative in \eqref{eq:ocp0b} using simple, local, first-order circular approximation.  The integral of a function $y$ is approximated using a trapezoidal quadrature rule, which is given by
\begin{equation}\label{eq:quad}
	\int_{0}^{T} y(t) dt \approx \sum_{k=1}^{N} y(t_k)w_k, \quad w_k=\frac{T}{N}
\end{equation}
for a circular domain. The derivative of a function $y$ is evaluated locally as a forward finite difference, with a time-periodic wrapping at the terminal time interval:
\begin{subequations}
\begin{align} 
    \frac{d}{dt} \wh{y}_N(t_k) \approx & \,\, (\bar{y}_{k+1}-\bar{y}_k)\cdot \frac{N}{T}, \,\,\, k=1,\ldots,N-1 \\
    \frac{d}{dt} &\wh{y}_N(t_N) \approx (\bar{y}_{1}-\bar{y}_{N})\cdot \frac{N}{T}.
\end{align}
\label{eq:deriv}
\end{subequations}
The derivative operator may be written in the form
\begin{equation}\label{eq:derivmat}
	\frac{d}{dt} \wh{y}_N(t_j) =\sum_{k=1}^N D_{jk}\bar{y}_k,
\end{equation}
a differentiation matrix $D$ has the entries $D_{kk}=-N/T$ for $k=1,\ldots,N$, $D_{k,k+1}=N/T$ for $k=1,\ldots,N-1$, $D_{N,1}=1$, and zero elsewhere.

Using \eqref{eq:Ix}, \eqref{eq:quad}, and \eqref{eq:derivmat}, the OCP \eqref{eq:ocp0} is transcribed as the following nonlinear program, in which the decision variables are vectors of the local function values $\bar{x}=(\bar{x}_1,\ldots,\bar{x}_N)$ and $\bar{u}=(\bar{u}_1,\ldots,\bar{u}_N)$:
\begin{subequations}
\begin{align}
\!\!\!\!\!	\min\ \ & \bar{J}(\bar{x},\bar{u})=\sum_{k=0}^N \cL(t_k,\bar{x}_k,\bar{u}_k)w_k \label{eq:ocp1a}\\
\!\!\!\!\!	{\rm s.t.}\ \ & \dS f\bp{t_i,\bar{x}_i,\sum_{k=1}^ND_{ik}\bar{x}_k,\bar{u}_i,\bar{h}_i}=0, & \!\!\!\!\!  i=1,\ldots,N \label{eq:ocp1b}\\
\!\!\!\!\!	& g(\bar{x}_k,\bar{u}_k,\bar{h}_k) \leq 0, & \!\!\!\!\!  k=1,\ldots,N \label{eq:ocp1d}
	\end{align}
\label{eq:ocp1}
\end{subequations}
We have eliminated the time-periodicity constraints \eqref{eq:ocp0d} in the formulation \eqref{eq:ocp1}, because they are implicit in the discretization scheme. It is possible to show that solutions to \eqref{eq:ocp1} converge to extrema of \eqref{eq:ocp0} as $N\to\infty$, using a similar approach as for pseudospectral schemes \cite{ruths11cdc}.

The scheme above does not provide an exact approximation in the case of certain polynomial functions, as can be shown for Legendre-Gauss schemes. However, the main advantage of the proposed ``circular'' time-discretization approach is sparsity, which reduces the number of terms in the constraint Jacobian by an order $\mathcal O (N)$, eliminates the need for additional constraints on the initial and terminal states, and creates a computationally well-posed problem in the case of time-periodic parameter functions $h$.  %By applying the above approximation approach, the nonlinear programming problems arising from problems of the form \eqref{prob:msw0} can be solved very rapidly, thereby enabling finer time-discretization and thus better fidelity of approximation.  In particular, terminal conditions \eqref{eq:termcon} can be omitted, thereby addressing significant numerical conditioning issues and ambiguity in defining initial and terminal conditions.  

%Practical implementation of optimization for gas pipeline operations, and in particular in a rolling horizon, model-predictive manner, will require assimilation of non-time-periodic boundary conditions.  
Real pipeline systems are subject to transient states over large spatiotemporal scales, any available baseline flow forecasts are uncertain and subject to modification, and measurements of system states are noisy.  In order to enable assimilation of data from the supervisory control and data acquisition (SCADA) system from a pipeline network into a model predictive OCP, we have tested our modeling in extensive validation studies \cite{zlotnik17psig}, and examined state estimation approaches \cite{sundar18tcst}.  To compensate for non-time-periodicity in future implementations of pipeline transient optimization and intra-day gas market mechanisms, we formulate a modified OCP over an extended time-horizon over which data are interpolated to produce periodic inputs, and where the solution can be taken as the restriction to the time-horizon of interest.  

Suppose now that the parameter functions $h$ in problem \eqref{eq:ocp0} are given as continuous on $[0,T]$, but not necessarily time-periodic. Let $\tau$ be an additional time constant by which the optimization horizon is extended, to be $[0,T+\tau]$.  We then construct parameter functions $\tilde{h}$ by linear interpolation on $[T,T+\tau]$, defined by $\tilde{h}(t)=h(t)$ for $t\in[0,T]$ and $\tilde{h}(t)=h(T)+(h(0)-h(T))(t-T)/\tau$ for $t\in[T,T+\tau]$.  Then a time-periodic problem of the form \eqref{eq:ocp0} can solved on the extended domain $[0,T+\tau]$, using the nonlinear programming formulation \eqref{eq:ocp1}. This results in time-periodic solutions $\tilde{x}$ and $\tilde{u}$ for the state and the control.  To obtain solutions for the state and control on the interval of interest, we simply use the restrictions $x(t)=\tilde{x}(t)$ and $u(t)=\tilde{u}(t)$ for $t\in[0,T]$.

Given a state solution that is consistent with boundary values, an instantaneous state can be used as an initial value constraint for a new time horizon that uses an update of these boundary conditions (in the future).  That is, once a solution is obtained given time-series for boundary data on a time interval $[T_1,T_1+T]$, a problem can be solved with updated time-series on a time interval $[T_1+H,T_1+H+T]$, where, e.g., $H\equiv1$ hour is the time interval over which the rolling horizon moves between solves.  The initial state for the latter problem can be constrained to the value obtained at $x(T_1+H)$ in the former problem, without losing the advantageous conditioning properties of a periodic formulation. 

\begin{figure}[t!]
\centering
\includegraphics[width=\linewidth]{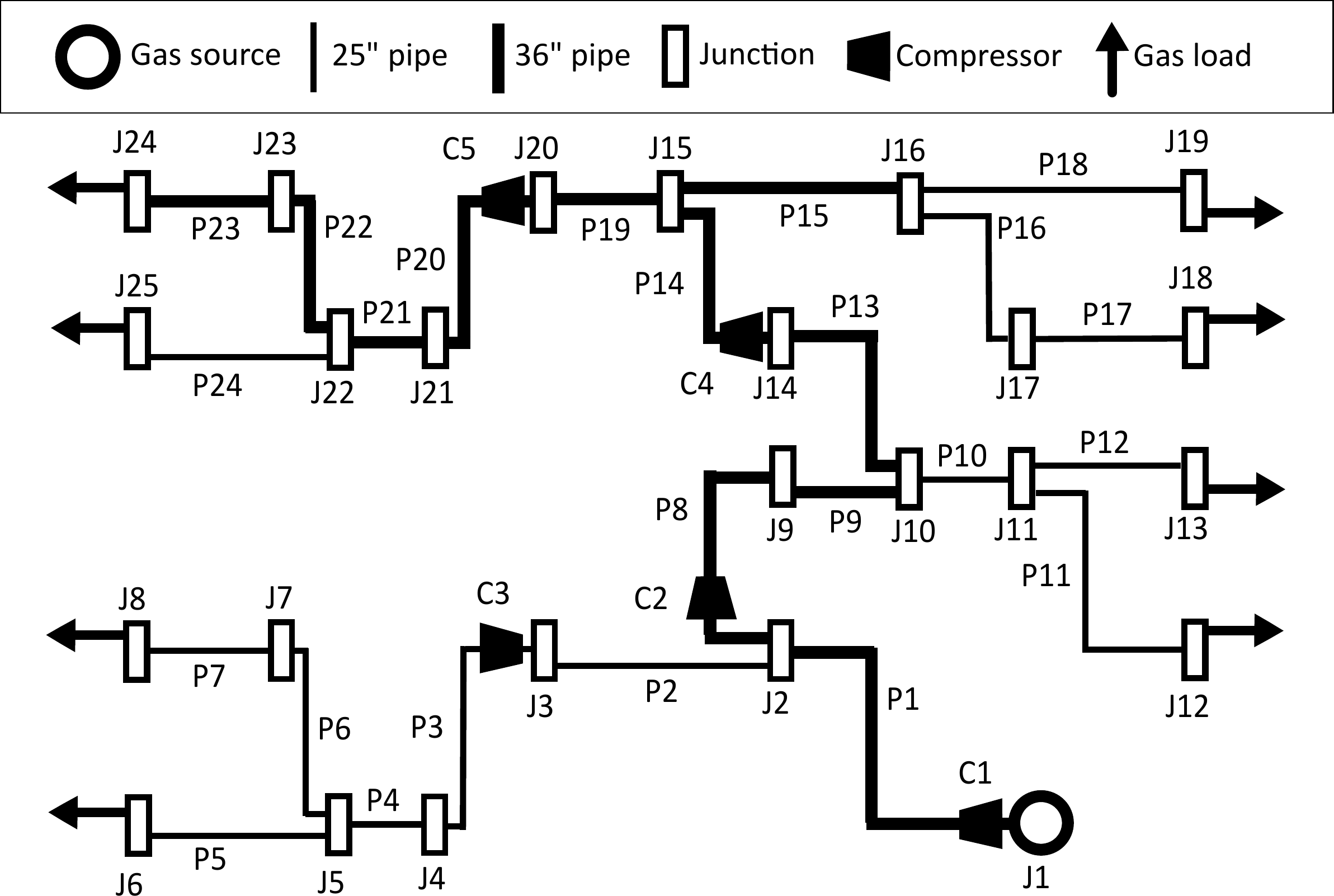}
\caption{LANL benchmark 25 node pipeline test network. Pipes (P1 to P24), pipe junctions (J1 to J25), and compressors (C1 to C5) are shown.}
\label{fig:gastestnetwork}
\end{figure}

\begin{figure}[t!]
\centering
\vspace{1ex}
\!\! \includegraphics[width=1.05\linewidth]{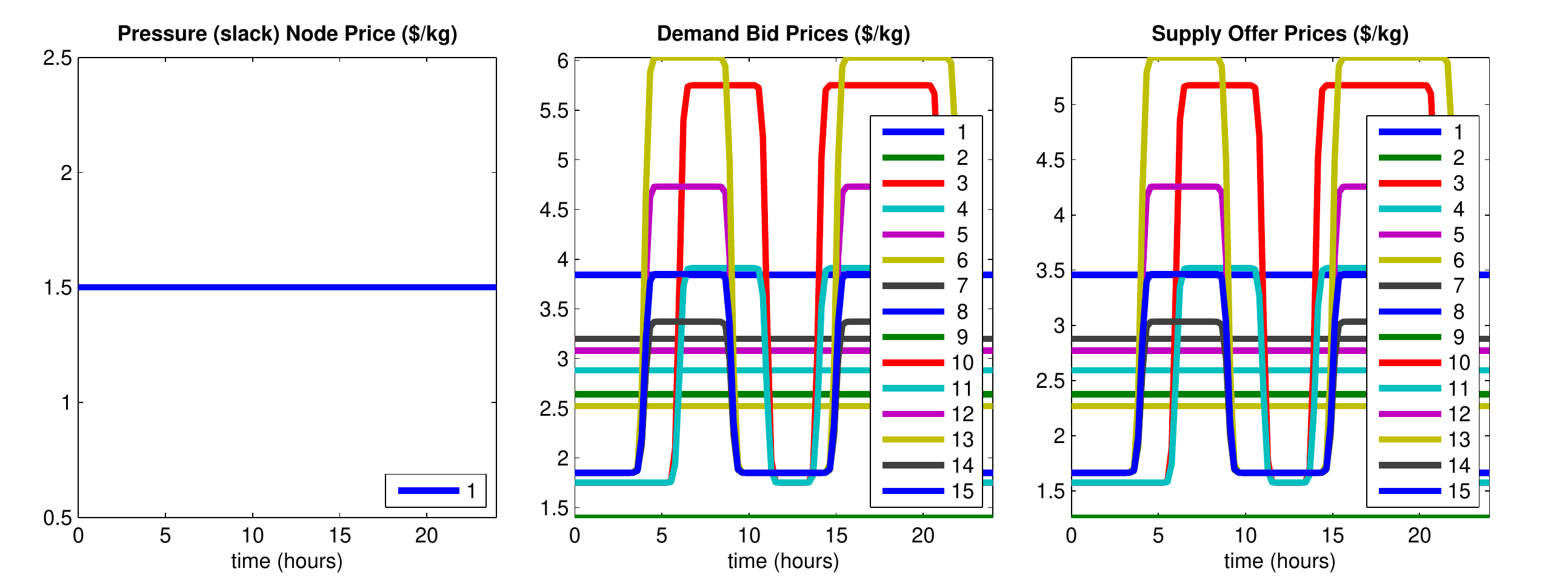}
\vspace{-2ex}
\!\! \includegraphics[width=1.05\linewidth]{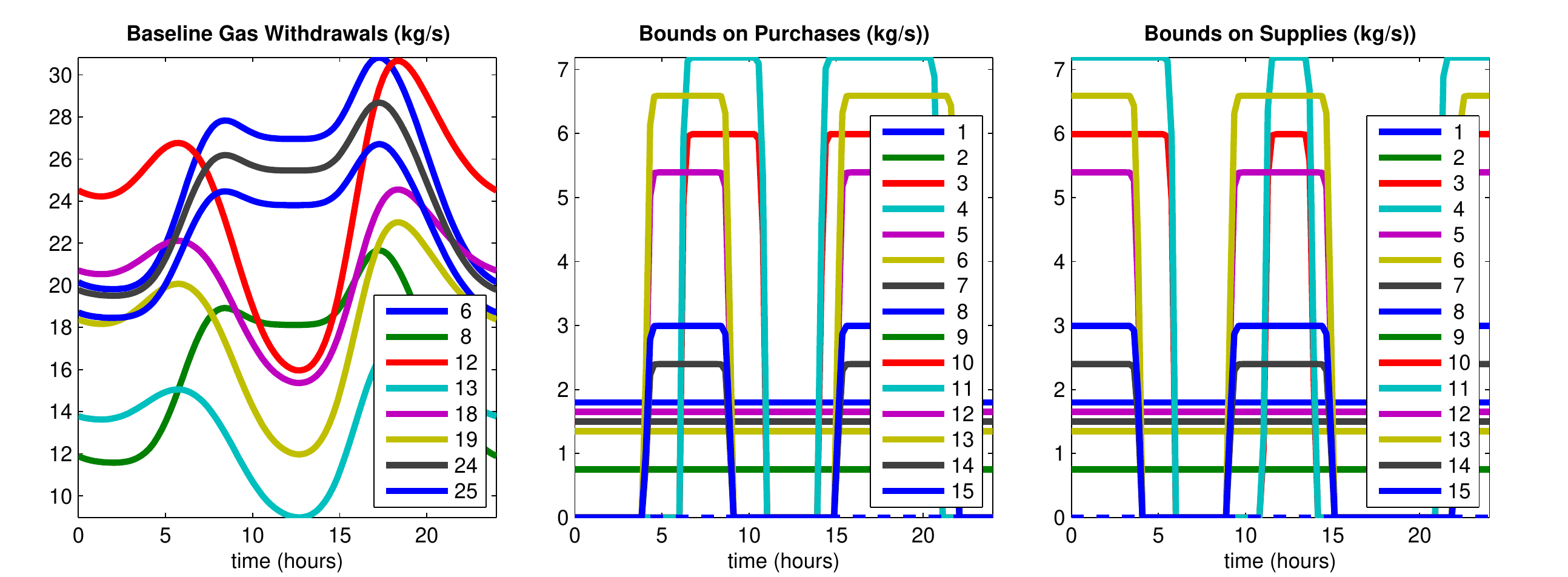}
\caption{Top: Price bids (Left: Constant price at slack node; Center: Buyer bids at transfer nodes;  Right: Seller offers at transfer nodes) Bottom: Quantity bids (Left: baseline flows; Center: Buyer demand bids; Right: Seller supply offers)}
\label{fig:c-bids}
\end{figure}

We apply the above optimization technique to formulate a rolling-horizon model predictive OCP for a gas pipeline auction market of the form \eqref{prob:msw0}.  We consider the internal system density $\bm \rho$ and flow $\bm q$ to describe the system state $x$, and the compressor ratios $\ubar{\alpha}_{ij}$ and $\bar{\alpha}_{ij}$ for $(i,j)\in\cC$ and the transfer node demands $d_m$ and supplies $s_m$ for $m\in\cG$ are the controls $u$.  The parameter functions $h$ include nominal nodal flows $\bar{q}_j$ for $j\in\cV_q$, supply density at slack nodes $\sigma_j$ for $j\in\cV_\sigma$, and prices $c_m^s(t)$ and $c_m^d(t)$ and constraint bound values $s_m^{\min}(t)$, $s_m^{\max}(t)$, $d_m^{\min}(t)$, and $d_m^{\max}(t)$ for $m\in\cG$.  The optimal control scheme is implemented by defining MATLAB functions for the objective, constraints, and their gradients with respect to decision variables, which are provided to the interior-point solver IPOPT version 3.11.8 running with the sparse linear solver ma57 \cite{biegler09ipopt}.  The method is available as an open source research code ``GRAIL'', which includes flexible routines for gas pipeline transient optimization and simulation \cite{zlotnik18grail}, for solving problems of operational or market designs and interfacing with power systems optimization software.  %The software features additional modeling details that are not described here, including non-ideal gas modeling using the CNGA formula for gas compressibility $Z$.   
The ``GRAIL'' tool has been used to evaluate the economic advantages of implementing an intra-day gas market, using a case study for a real pipeline system and associated SCADA time-series data \cite{rudkevich19hicss}.  Key inputs and outputs are examined in the case study.

\begin{figure}[t!]
\centering
\!\! \includegraphics[width=1.05\linewidth]{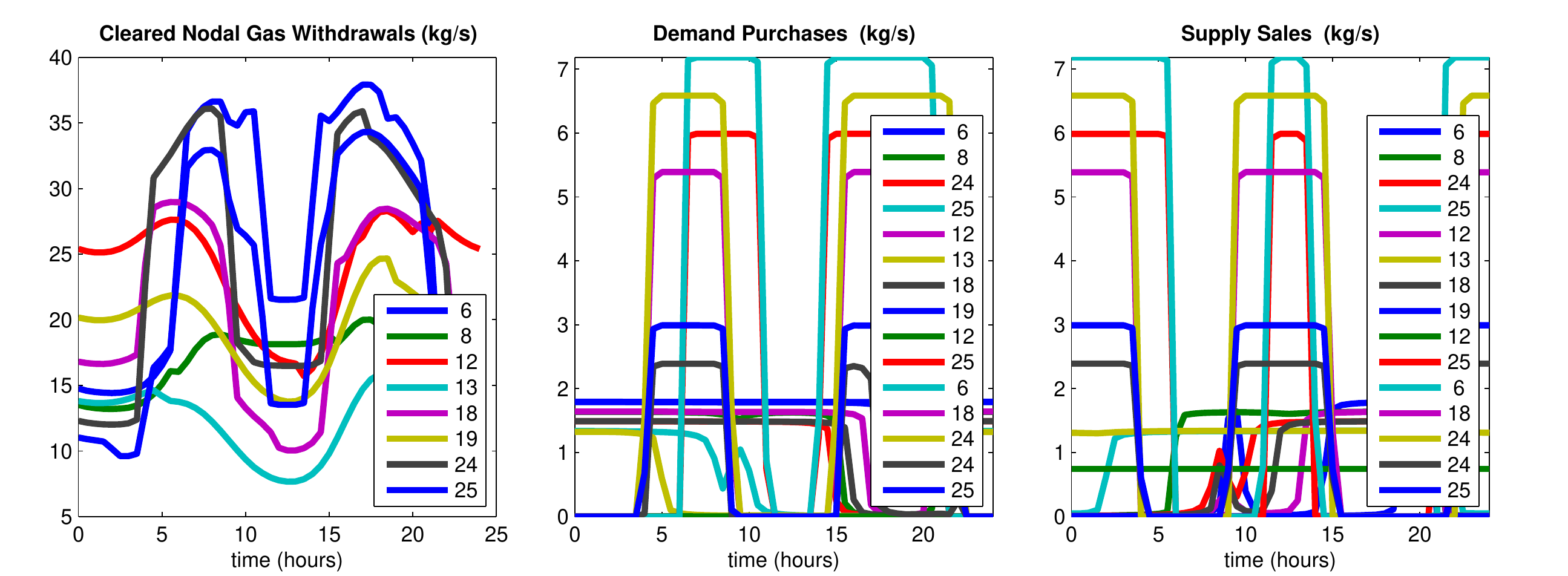}
\vspace{-2ex}
\!\! \includegraphics[width=1.05\linewidth]{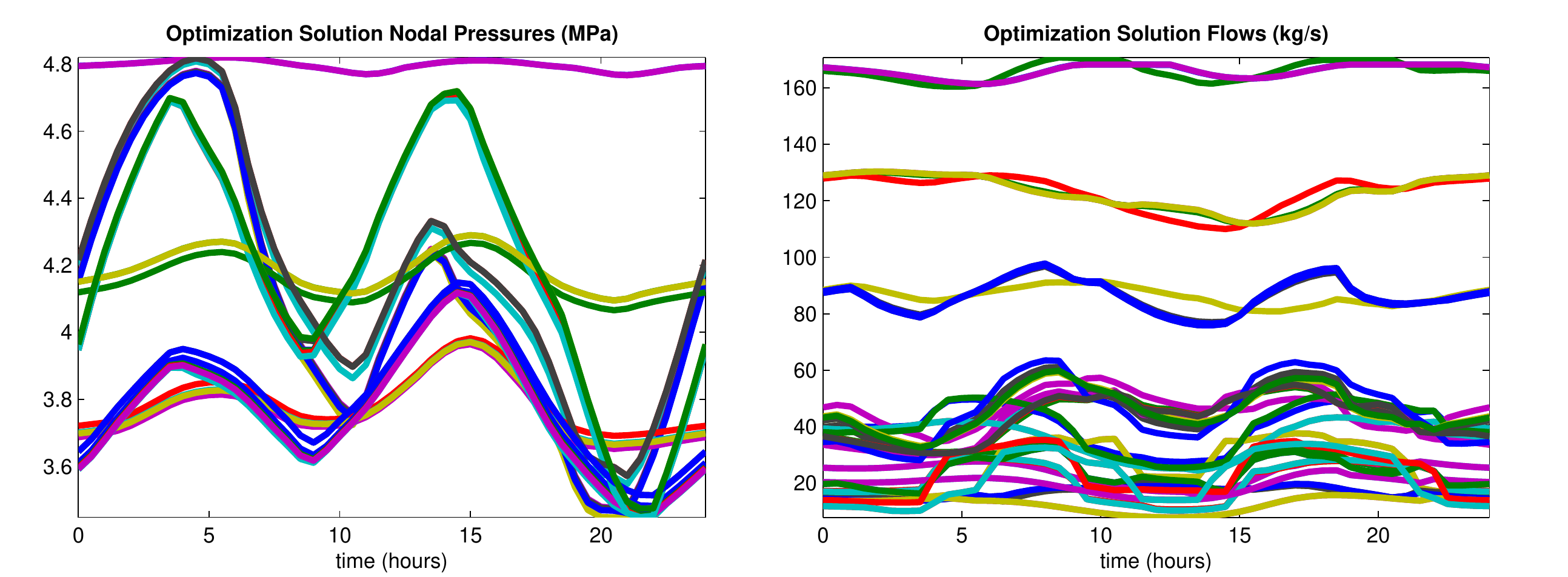}
\vspace{-2ex}
\!\! \includegraphics[width=1.05\linewidth]{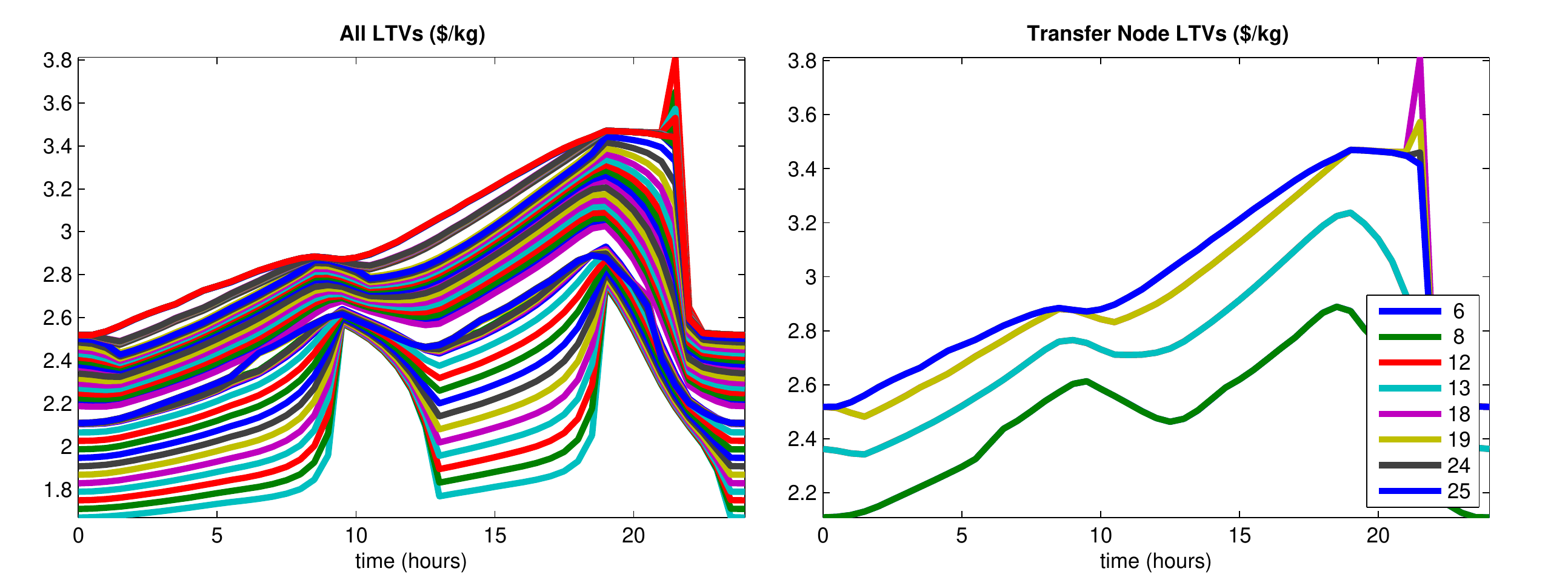}
\caption{Top row: Flow schedules (Left: total flows to physical locations with transfer nodes; Center: purchases at each transfer node; Right: sales at each transfer node). Middle row: Internal system state (Left: Pressure at physical nodes, flows at pipe inlets and outlets). Bottom row: Price schedule (Left: prices at all system nodes; Right: prices at transfer nodes).}
\label{fig:ltv}
\end{figure}

\section{Case Study} \label{sec:example}

We present a (time-periodic) case study for clearing an intra-day market for a standard pipeline test network (see Figure \ref{fig:gastestnetwork}), which was  used in previous studies \cite{zlotnik15cdc,mak16acc}.  Notably, the low order, local scheme maximizes sparsity of the NLP constraint Jacobian.  Here the Jacobian of the NLP in the case study has under 0.0745\% non-zero entries, and the solution requires less than 20 seconds on a commodity computer using 10 km spatial discretization and 24 collocation points over 24 hours.  While the entire model inputs and outputs cannot be fully presented here, we encourage the reader to view the full case study results available as an example with the GRAIL software \cite{zlotnik18grail}.  The market bids are shown in Figure \ref{fig:c-bids}, and physical and price solutions are shown in Figure \ref{fig:ltv}.  The key observation is that the demand bids shown in Fig. \ref{fig:c-bids} top center cannot be entirely fulfilled (see Fig. \ref{fig:ltv} top center), and binding constraints result in price separation (Fig. \ref{fig:ltv} bottom right) as in the steady-state \cite{rudkevich17hicss}.

\section{Conclusion} \label{sec:conc}

We have presented an economic optimal control problem for scheduling and pricing natural gas flows in pipeline transmission systems, as well as a highly efficient approximation scheme solving the problem for realistic systems.  The method was applied the to a pipeline system test network case study, in which time-dependent flow schedules and locational trade values of gas were computed. Future work will involve extension to mixed-integer formulations \cite{gugat18mip}, incorporate recent modeling advances \cite{hante17}, and transition to practice for real systems \cite{zlotnik19psig}.

\section*{Acknowledgement} This work was carried out as part of Project GECO for the Advanced Research Project Agency-Energy of the U.S. Department of Energy under Award No. DE-AR0000673. Work at Los Alamos National Laboratory was conducted for the D.O.E. Office of Electricity Advanced Grid Research and Development program under the auspices of the National Nuclear Security Administration of the U.S. Department of Energy under Contract No. 89233218CNA000001.

\bibliographystyle{unsrt}
\bibliography{gas_master,markets_master,spe,power_master}

\end{document}

%% file: avz_defs.tex
\usepackage{amsmath, amssymb, amsxtra, mathrsfs, bm}
\usepackage{graphics, epsfig, here, epstopdf, fancyhdr,graphicx}
\usepackage{bbm}
\usepackage{color}
\usepackage{multirow}
\usepackage[all]{xy}
\usepackage{accents}

\newcommand{\ubar}[1]{\underaccent{\bar}{#1}}

\newcommand{\bR}{\mathbb{R}}

\newcommand{\cC}{\mathcal{C}}
\newcommand{\cE}{\mathcal{E}}

\newcommand{\cL}{\mathcal{L}}

\newcommand{\cV}{\mathcal{V}}

\newcommand{\cT}{\mathcal{T}}

\newcommand{\fA}{\,\forall\,}

\newcommand{\dS}{\displaystyle}

\newcommand{\bp}[1]{{\left(#1\right)}}

\newcommand{\teQ}{\triangleq}

\newcommand{\rM}[1]{{\mathrm{#1}}}

\newcommand{\cG}{\mathcal{G}}

\newcommand{\wh}[1]{{\widehat{#1}}}

\newcommand{\rd}{\mathrm{d}}